\RequirePackage{stix}
\documentclass[11pt, reqno]{amsart}

\usepackage{amsthm, amsmath, amsfonts}
\usepackage[alphabetic]{amsrefs}
\usepackage{verbatim}
\usepackage{graphicx}

\newcommand{\newword}{\textbf}

\newcommand{\demo}{\begin{proof}}
\newcommand{\edemo}{\end{proof}}
\newcommand{\demoname}[1]{\begin{proof}[#1]}
\newcommand{\edemoname}{\end{proof}}

\newcommand{\stepname}{Step}

\theoremstyle{plain}
\newtheorem{theorem}{Theorem}[section]

\newtheorem{conjecture}[theorem]{Conjecture}
\newtheorem{conjlemma}[theorem]{Conjectured Lemma}
\newtheorem{corollary}[theorem]{Corollary}
\newtheorem{lemma}[theorem]{Lemma}

\theoremstyle{definition}
\newtheorem{example}[theorem]{Example}
\newtheorem{definition}[theorem]{Definition}
\newtheorem{remark}[theorem]{Remark}
\newtheorem{conjremark}[theorem]{Conjectured Remark}

\newcommand{\thm}{\begin{theorem}}
\newcommand{\ethm}{\end{theorem}}
\newcommand{\conj}{\begin{conjecture}}
\newcommand{\econj}{\end{conjecture}}
\newcommand{\conjlem}{\begin{conjlemma}}
\newcommand{\econjlem}{\end{conjlemma}}
\newcommand{\conjremk}{\begin{conjremark}}
\newcommand{\econjremk}{\end{conjremark}}

\newcommand{\expl}{\begin{example}}
\newcommand{\eexpl}{\qex\end{example}}
\newcommand{\defn}{\begin{definition}}
\newcommand{\edefn}{\qef\end{definition}}
\newcommand{\defnnb}{\begin{definition}}
\newcommand{\edefnnb}{\end{definition}}
\newcommand{\remk}{\begin{remark}}
\newcommand{\eremk}{\qex\end{remark}}
\newcommand{\remknb}{\begin{remark}}
\newcommand{\eremknb}{\end{remark}}

\newcommand{\coro}{\begin{corollary}}
\newcommand{\ecoro}{\end{corollary}}
\newcommand{\lem}{\begin{lemma}}
\newcommand{\elem}{\end{lemma}}

\providecommand{\qexsymbol}{$\lozenge$}%
\newcommand{\mathqex}{\quad\hbox{\qexsymbol}}
\DeclareRobustCommand{\qex}{%
  \ifmmode \mathqex
  \else
    \leavevmode\unskip\penalty9999 \hbox{}\nobreak\hfill
    \quad\hbox{\qexsymbol}%
  \fi
}
\providecommand{\qefsymbol}{$\triangle$}%
\newcommand{\mathqef}{\quad\hbox{\qefsymbol}}
\DeclareRobustCommand{\qef}{%
  \ifmmode \mathqef
  \else
    \leavevmode\unskip\penalty9999 \hbox{}\nobreak\hfill
    \quad\hbox{\qefsymbol}%
  \fi
}

\newcommand{\enum}{\begin{enumerate}}
\newcommand{\eenum}{\end{enumerate}}

\newcounter{countalph}
\newenvironment{alphnumerate}%
{\begin{list}{{\rm (\alph{countalph})}}{\usecounter{countalph}}}%
{\end{list}}
\newcommand{\enumalph}{\begin{alphnumerate}}
\newcommand{\eenumalph}{\end{alphnumerate}}

\newcommand{\rc}{\mathop{\circ}}

\newcommand{\rrr}[1]{{\mathbb{R}}^{#1}}

\newcommand{\rw}{\rrr{2}}
\newcommand{\rt}{\rrr{3}}
\newcommand{\zz}{\mathbb{Z}}

\newcommand{\nn}{\mathbb{N}}

\newcommand{\ztlat}{${\mathbb{Z}}^3$ lattice}
\newcommand{\zwlat}{${\mathbb{Z}}^2$ lattice}

\newcommand{\xyplane}{\nd{xy}plane}

\newcommand{\lif}{lift}
\newcommand{\lifs}{lifts}
\newcommand{\xlif}{\nd{x}lift}
\newcommand{\xlifs}{\nd{x}lifts}
\newcommand{\ylif}{\nd{y}lift}

\newcommand{\xlift}[1]{\nd{(x, #1)}lift}

\newcommand{\ylift}[1]{\nd{(y, #1)}lift}

\newcommand{\prlif}{proper}
\newcommand{\plif}{\prlif\ lift}
\newcommand{\plifs}{\prlif\ lifts}
\newcommand{\pxlif}{\prlif\ \nd{x}lift}
\newcommand{\pxlifs}{\prlif\ \nd{x}lifts}
\newcommand{\pylif}{\prlif\ \nd{y}lift}
\newcommand{\pylifs}{\prlif\ \nd{y}lifts}

\newcommand{\pylift}[1]{\prlif\ \nd{(y, #1)}lift}

\newcommand{\func}[3]{{#1} \colon {#2} \to {#3}}

\newcommand{\nd}[1]{$#1$\nobreakdash-\hspace{0pt}}
\newcommand{\nb}[1]{{#1}\nobreakdash-\hspace{0pt}}

\newcommand{\simpoab}[2]{\langle {#1}, {#2}\rangle}

\newcommand{\prop}{regular}
\newcommand{\lsprop}{proper}

\newcommand{\ralp}{lattice diagram of a knot or link}

\newcommand{\ralps}{lattice diagrams of knots and links}

\newcommand{\ralpk}{lattice diagram of a knot}

\newcommand{\ralpl}{lattice diagram of a link}

\newcommand{\ral}{lattice diagram}

\newcommand{\lsk}{lattice stick knot}

\newcommand{\lskl}{lattice stick knot or link}

\newcommand{\lskls}{lattice stick knots and links}
\newcommand{\Lskls}{Lattice stick knots and links}
\newcommand{\LSKLs}{Lattice Stick Knots and Links}

\newcommand{\lsplsk}{\lsprop\ lattice stick knot}

\newcommand{\lsplskl}{\lsprop\ lattice stick knot or link}

\newcommand{\lsplskls}{\lsprop\ lattice stick knots and links}

\newcommand{\aplskl}{lattice pre-image}

\newcommand{\zvalue}{\nd{z}value}
\newcommand{\zvalues}{\nd{z}values}

\newcommand{\heig}{\zvalue}

\newcommand{\zstick}{\nd{z}stick}
\newcommand{\zsticks}{\nd{z}sticks}

\newcommand{\xzstick}{\nd{x}facing \nd{z}stick}
\newcommand{\xzsticks}{\nd{x}facing \nd{z}sticks}
\newcommand{\yzstick}{\nd{y}facing \nd{z}stick}
\newcommand{\yzsticks}{\nd{y}facing \nd{z}sticks}
\newcommand{\uzstick}{upper \nd{z}stick}

\newcommand{\lzstick}{lower \nd{z}stick}

\newcommand{\xstrand}{\nd{x}strand}

\newcommand{\ystrand}{\nd{y}strand}

\newcommand{\xpstrand}{\nd{x}pre-strand}

\newcommand{\ypstrand}{\nd{y}pre-strand}

\newcommand{\xarc}{\nd{x}arc}
\newcommand{\xarcs}{\nd{x}arcs}
\newcommand{\yarc}{\nd{y}arc}
\newcommand{\yarcs}{\nd{y}arcs}
\newcommand{\zarc}{\nd{z}arc}
\newcommand{\zarcs}{\nd{z}arcs}

\newcommand{\xzarc}{\nd{x}facing \nd{z}arc}
\newcommand{\yzarc}{\nd{y}facing \nd{z}arc}

\newcommand{\xedge}{\nd{x}edge}
\newcommand{\xedges}{\nd{x}edges}
\newcommand{\yedge}{\nd{y}edge}
\newcommand{\yedges}{\nd{y}edges}
\newcommand{\zedge}{\nd{z}edge}
\newcommand{\zedges}{\nd{z}edges}

\newcommand{\zproj}{\pi_{xy}}

\newcommand{\nearby}{nearby}
\newcommand{\oppsite}{opposing}

\newcommand{\vcrossing}{\nd{y}crossing}
\newcommand{\vcrossings}{\nd{y}crossings}
\newcommand{\hcrossing}{\nd{x}crossing}
\newcommand{\hcrossings}{\nd{x}crossings}

\newcommand{\slpc}{Celtic configuration}

\newcommand{\probcrs}{problem crossing}
\newcommand{\probcrss}{problem crossings}\newcommand{\nprobcrs}{non-problem crossing}
\newcommand{\badnbr}{bad neighbor}
\newcommand{\badnbrs}{bad neighbors}
\newcommand{\crnr}{corner}

\newcommand{\latgr}{lattice graph}

\newcommand{\crosgr}{crossing graph}

\newcommand{\CROSGR}{Crossing Graph}

\newcommand{\pcrosgr}{problem crossing graph}

\newcommand{\PCROSGR}{Problem Crossing Graph}

\newcommand{\crg}[1]{CG_{#1}}

\newcommand{\pcrg}[1]{PG_{#1}}

\newcommand{\usf}{deleted-square free}

\newcommand{\latpath}{lattice arc}
\newcommand{\latpaths}{lattice arcs}
\newcommand{\latcyc}{lattice simple closed curve}
\newcommand{\latcycs}{lattice simple closed curves}

\newcommand{\cmpt}{compatible}
\newcommand{\bcmpt}{backwards compatible}

\newcommand{\latmap}{lattice map}

\newcommand{\resv}{resolved}
\newcommand{\nresv}{\nb{non}resolved}

\newcommand{\pinvx}[1]{\zproj^{-1}({#1})}

\newcommand{\twoless}{\nd{2}near regular}

\newcommand{\tworeg}{\nd{2}regular}

\newcommand{\thmref}[1]{Theorem~\ref{#1}}
\newcommand{\lemref}[1]{Lemma~\ref{#1}}

\newcommand{\defnref}[1]{Definition~\ref{#1}}

\newcommand{\lemrefp}[2]{Lemma~\ref{#1}~(\ref{#1#2})}

\newcommand{\defnrefp}[2]{Definition~\ref{#1}~(\ref{#1#2})}

\newcommand{\figref}[1]{Figure~\ref{#1}}
\newcommand{\figrefs}[2]{Figures~\ref{#1} and \ref{#2}}

\newcommand{\remkrefp}[2]{Remark~\ref{#1}~(\ref{#1#2})}

\begin{document}

\title{Lattice Diagrams of Knots and Diagrams of Lattice Stick Knots}
\author{Margaret Allardice}
\author{Ethan D.\ Bloch}
\address{Bard College\\
Annandale-on-Hudson, NY 12504\\
U.S.A.}
\email{bloch@bard.edu}
\begin{abstract}
We give a simple example showing that a knot or link diagram that lies in the \zwlat\ is not necessarily the projection of a \lskl\ in the \ztlat, and we give a necessary and sufficient condition for when a knot or link diagram that lies in the \zwlat\ is in fact the projection of a \lskl.
\end{abstract}
\maketitle

\section{Introduction}

\Lskls, that is, knots and links that are in the \ztlat\ (which is the graph in $\rt$ where the vertices are the points with integer coefficients, and the edges are unit length and parallel to the coordinate axes), have been studied by a number of authors, for example \cite{DIAO1}, \cite{ER-PH}, \cite{JVR-P}, \cite{U-J-O-T-W}, \cite{S-I-A-D-S-V}, \cite{H-H-K-N-O}, \cite{D-E-P-Z}, \cite{D-E-Y} and \cite{H-K-O-N}.  There is some variation in terminology in these and other papers; for example, some authors use the term ``cubic lattice'' rather than \ztlat, and some use ``step'' to mean an edge in the \ztlat.  The general goal is to find the minimum number of edges needed to represent a given knot or link as a \lskl, and, when the minimum number of edges has not been found, to give an upper bound for it.

Here we address a more basic question, inspired by an analogous comment in the stick knot (but not \lsk) case in \cite{AD-SH}, which defines the ``projection stick index'' of a knot $K$ as ``the least number of sticks in any projection of a polygonal conformation of $K$,'' and where they add ``Note that there is no \textit{a priori} reason that it is equal to the least number of sticks in polygonal projections of knot embeddings that are not themselves polygonal.''

An analogous question could be asked in the case of stick numbers and edge numbers of projections of \lskls, but, more fundamentally, we raise the question of the relation between projections of \lskls\ on the one hand, and, on the other hand, projections of knots and links in $\rt$ that are not necessarily themselves \lskls\ but where the projection is in the \zwlat.

For standard knots and links, there is no analogous question to be asked.  That is, any knot or link diagram in $\rw$ is the \prop\ projection of a knot or link in $\rt$.  That implies, for example, that if a knot or link in $\rt$ is projected onto $\rw$, and if the resulting diagram is then manipulated (e.g.\ using the Reidemeister moves), then after the manipulation the diagram is still the projection of a knot or link in $\rt$.

The situation is not the same in the case of \lskls.  In particular, we give a simple example showing that of a knot or link diagram that lies in the \zwlat\ is not necessarily the projection of a \lskl\ in the \ztlat.  In \thmref{thmBDM}, we give a necessary and sufficient condition for when a knot or link diagram that lies in the \zwlat\ is in fact the projection of a \lskl.


\section{Projections of \LSKLs}

We start with some terminology.  An \newword{\xedge}, \newword{\yedge} and \newword{\zedge} of the \zwlat\ or \ztlat\ is an edge that is parallel to the corresponding coordinate axes.  A \newword{\zstick} in the \ztlat\ is a line segment that is the union of finitely many \zedges; for convenience, we will sometimes consider vertices in the \zwlat\ to be \newword{trivial} \zsticks.

Let $\zproj$ denote orthogonal projection from $\rt$ onto the \xyplane.

\begin{figure}[ht]
\centering\includegraphics[scale=0.3]{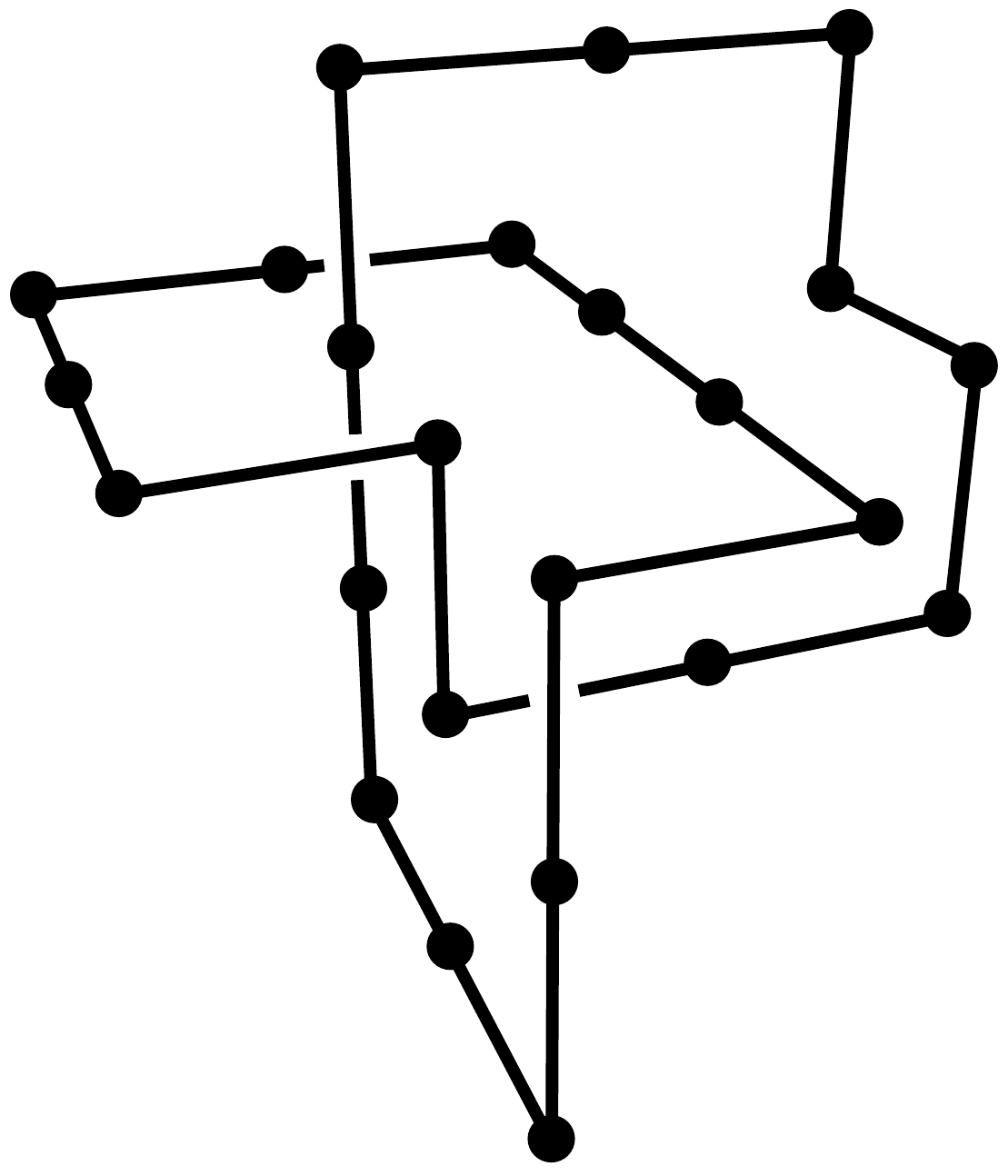}
\caption{}\label{ralp_knot_3_1_3f_rev}
\end{figure}

The objects we wish to study are \newword{\lskls}, which are knots and links contained in the \ztlat.  See \figref{ralp_knot_3_1_3f_rev} for a \lsk\ representing the trefoil knot.

Recall that when forming the diagram of a knot or link, it is required that the knot or link be positioned in $\rt$ so that its projection onto the \xyplane\ is a \prop\ projection, which means that the inverse image under projection onto the \xyplane\ of any point in the diagram that is not a crossing is a single point, and the inverse image at a crossing (of which there are only finitely many) is two points.

\begin{figure}[ht]
\centering\includegraphics[scale=0.4]{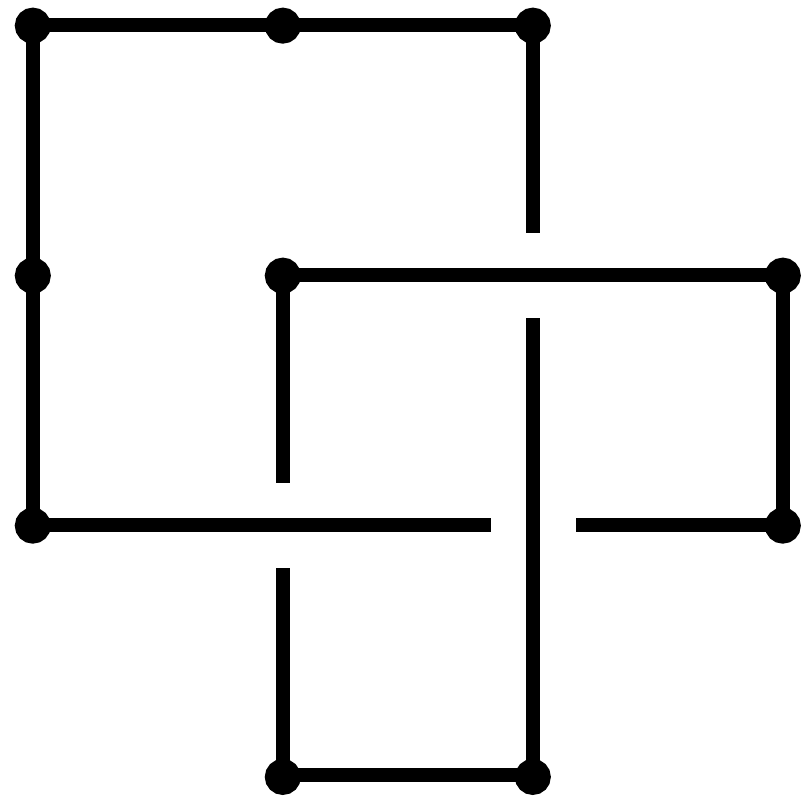}
\caption{}\label{ralp_knot_3_1_3i_rev}
\end{figure}

On the other hand, whereas the projection of a \lskl\ onto the \xyplane\ is indeed a diagram of a knot of link (see, for example, the \lsk\ in \figref{ralp_knot_3_1_3f_rev} and its projection in \figref{ralp_knot_3_1_3i_rev}), such a projection is not a \prop\ projection of a knot or link, because the inverse image of a vertex in the projection that is not crossing can be a non-trivial \zstick\ rather than a single point, and the the inverse image of a vertex in the projection that is a crossing can be one or two non-trivial \zsticks\ rather than two points.

That said, just as a \prop\ projection of a knot or link has restrictions on the possible inverse images of points, so too when we we project a \lskl\ onto the \xyplane.

\defnnb
Let $K$ be a \lskl.  
\enum
\item
The \lskl\ $K$ is \newword{\lsprop} if the inverse image under projection onto the \xyplane\ of any point in the diagram that is not a vertex of the \zwlat\ is a single point; the inverse image of a vertex in the diagram that is not a crossing is a single \zstick\ (possibly trivial); and the inverse image of a vertex in the diagram that is at a crossing is two \zsticks\ (again, possibly trivial).
\item
Suppose that $K$ is \lsprop.  Let $R$ be the projection of $K$, and let $a$ be a vertex of $R$.  Then $\pinvx a$ is the union of two \zsticks, one of which has larger \zvalues\ for all its points than the \zvalues\ of the other \zstick; the former of the two \zsticks\ is called the \newword{\uzstick} of $K$ over $a$, and the other \zstick\ is called the \newword{\lzstick}.
\qef
\eenum
\edefnnb

\begin{figure}[ht]
\centering\includegraphics[scale=0.4]{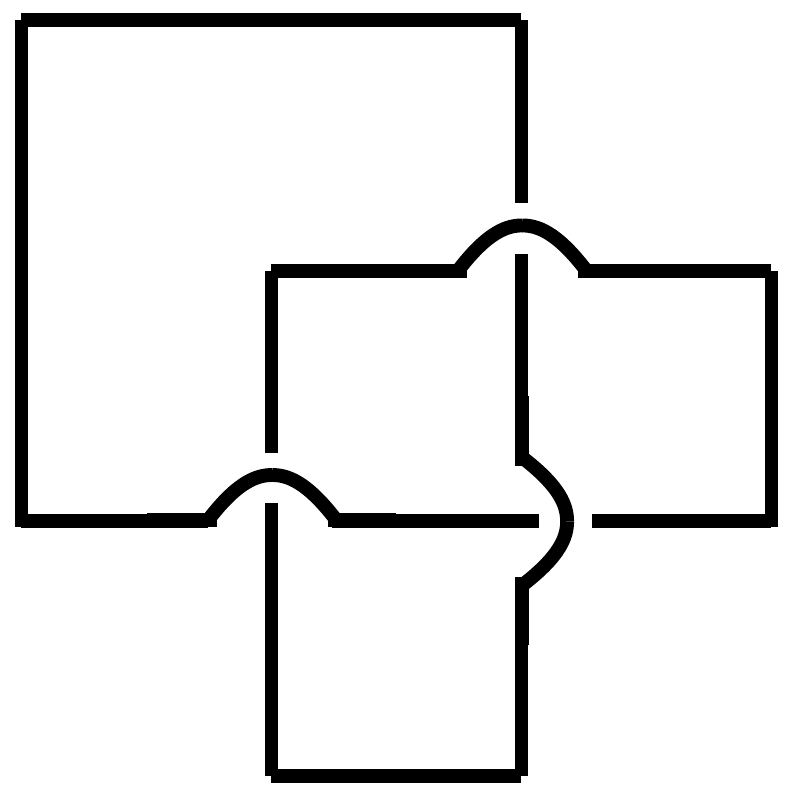}
\caption{}\label{ralp_knot_3_1_3k_rev}
\end{figure}

We note that whereas the knot diagram seen in \figref{ralp_knot_3_1_3i_rev} is naturally thought of as created by a (non-\prop) projection of a \lsplskl, this knot diagram can also be obtained by a \prop\ projection of a non-lattice knot, as seen in \figref{ralp_knot_3_1_3k_rev}.

The question we address is the reverse of the above observation, which we state using the following terminology.

\defn
A \newword{\ralp} is a diagram of a knot or link, such that the diagram is contained the the \zwlat, and that all crossings in the knot or link diagram are at vertices of the \zwlat.
\edefn

For example, the knot diagram in \figref{ralp_knot_3_1_3i_rev} is a \ralpk.

We then ask, is every \ralp\ the projection of a \lsplskl?  It might be thought that the answer is trivially yes, because we could simply do a lattice analog of what we did in the non-lattice context in \figref{ralp_knot_3_1_3k_rev}.  However, the following example shows that that is not always possible.

\begin{figure}[ht]
\centering\includegraphics[scale=0.4]{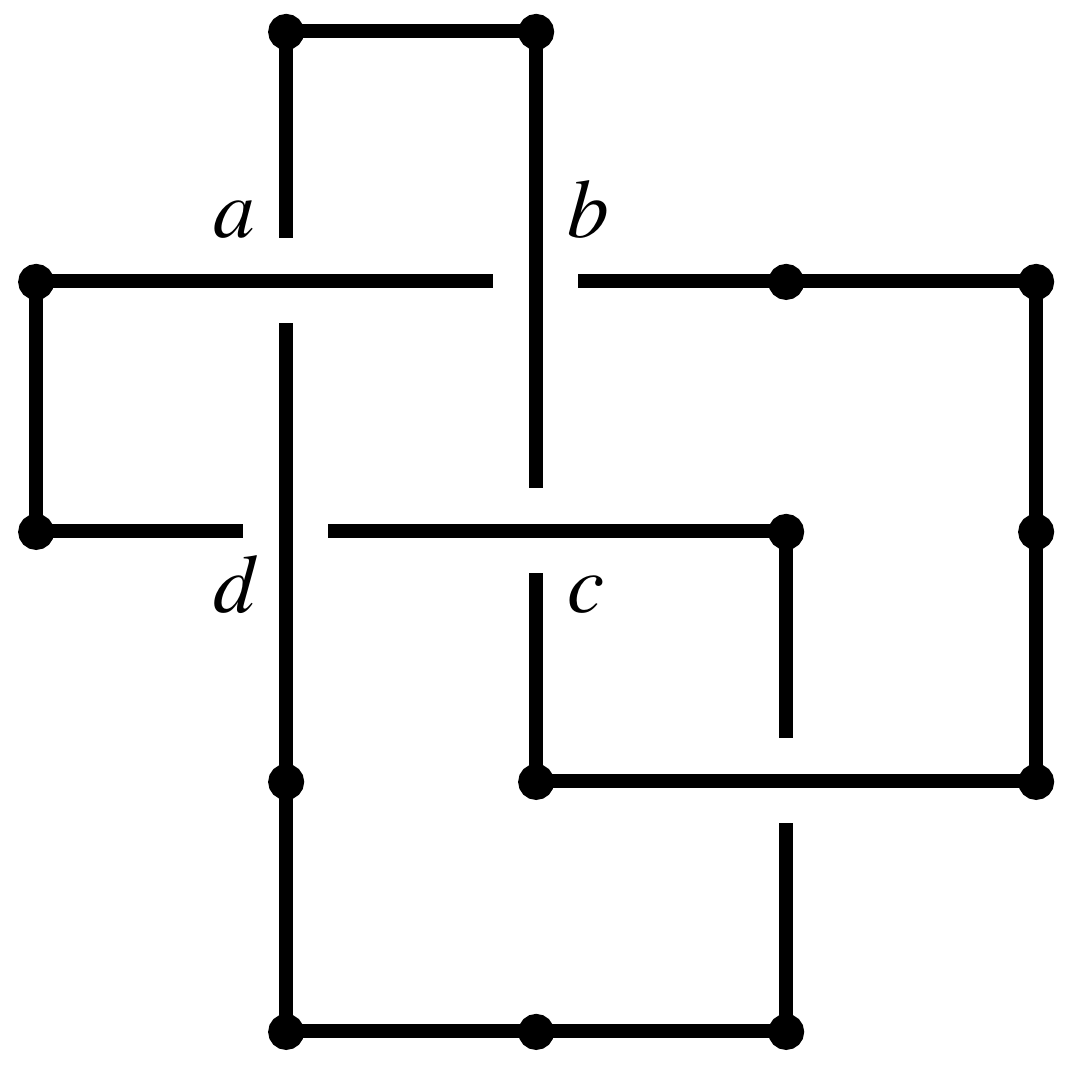}
\caption{}\label{ralp_knot_5_2a}
\end{figure}

Consider the \ral\ $R$ in \figref{ralp_knot_5_2a}, which is the knot $5_2$.  Suppose that $R$ is the projection of a \lsk\ $K$.  The edge of $K$ that projects onto $\simpoab ab$ is connected to the \lzstick\ of $K$ over $b$, and the edge of $K$ that projects onto $\simpoab bc$ is connected to the \uzstick\ of $K$ over $b$.  Hence the edge of $K$ that projects onto $\simpoab ab$ has smaller \zvalue\ than the edge of $K$ that projects onto $\simpoab bc$.  The same argument shows that the edge of $K$ that projects onto $\simpoab bc$ has smaller \zvalue\ than the edge of $K$ that projects onto $\simpoab cd$, that the edge of $K$ that projects onto $\simpoab cd$ has smaller \zvalue\ than the edge of $K$ that projects onto $\simpoab da$, and that the edge of $K$ that projects onto $\simpoab da$ has smaller \zvalue\ than the edge of $K$ that projects onto $\simpoab ab$, which leads to the obvious contradiction.  Hence $R$ is not the projection of a \lsplsk\ $K$.

Clearly, any \ralp\ with the same type of configuration as the four vertices $a$, $b$, $c$ and $d$ in \figref{ralp_knot_5_2a} is not be a projection of a \lsplskl.  As we will see in Theorem~\ref{thmBDM} below, such a configuration is the only obstacle to a \ralp\ being a projection of a \lsplskl.

\begin{figure}[ht]
\centering\includegraphics[scale=0.3]{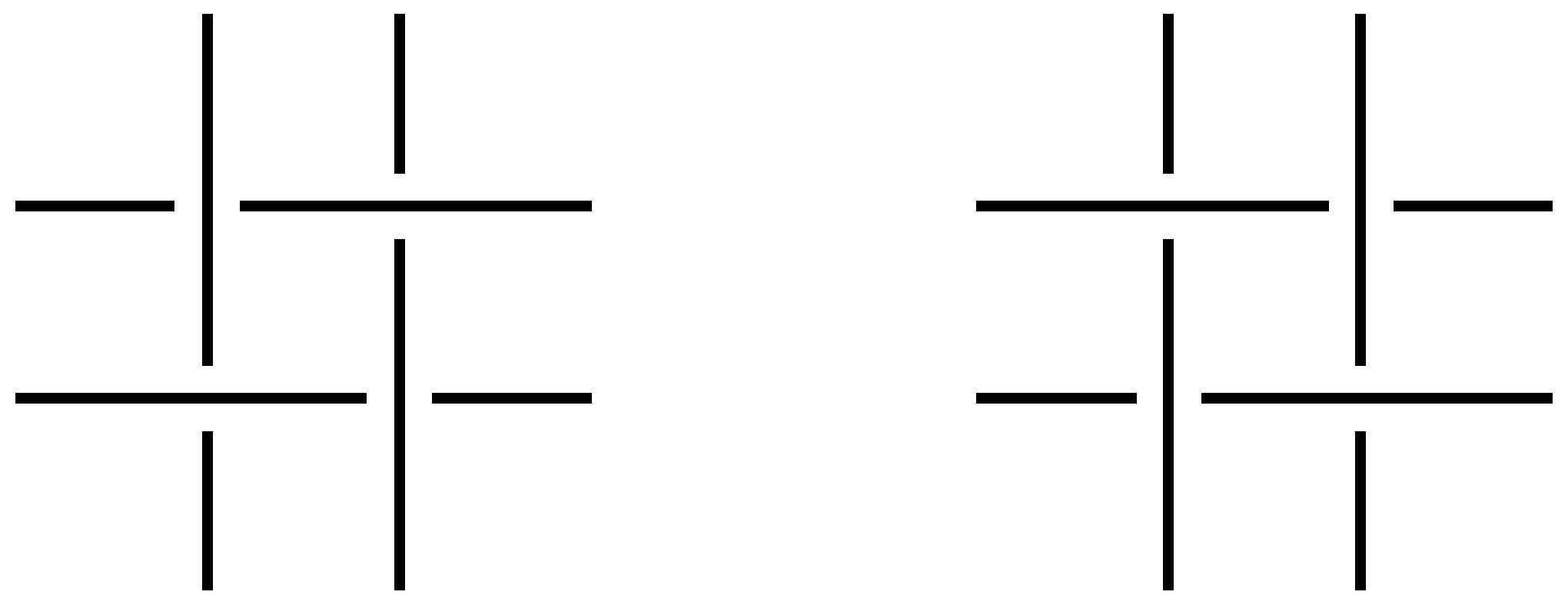}
\caption{}\label{ralp_knot_celtic}
\end{figure}

\defn\label{defnBDD}
A \newword{\slpc} is a subset of a \ralp\ that is equivalent to either of the configurations in \figref{ralp_knot_celtic}, where the crossings in the figure are at adjacent vertices in the \zwlat.
\edefn

The name ``\slpc'' is due to the fact that this configuration occurs regularly (albeit often on the diagonal) in Celtic interlace patterns, as in \cite{FISH}, for example; it also occurs in interlace patterns in other cultures, for example China.

Our theorem is the following.

\thm\label{thmBDM}
Let $R$ be a \ralp.  Then $R$ is the projection of a \lsplskl\ if and only if $R$ does not have a \slpc.
\ethm

The proof of Theorem~\ref{thmBDM} will be given in Section~\ref{secBASI}, after some preliminaries.


\section{\CROSGR\ and \PCROSGR}

We now define two graphs that arise from \ralps, using the following terminology.  A \newword{\latgr} is a subgraph the \zwlat; a \newword{\latpath}, respectively \newword{\latcyc}, is a \latgr\ that is a path graph, respectively cycle graph.

Note that if $G$ is a \latgr, and if $a$ and $b$ are vertices of $G$, we say that $a$ and $b$ are ``adjacent'' if they are joined by an edge of $G$, not if they are joined only by an edge in the \zwlat.

\defnnb\label{defnBEE}
Let $G$ be a \latgr, and let $v$ be a vertex of $G$.
\enum
\item
The \latgr\ $G$ is \newword{\usf} if it does not have three edges that are part of a unit square in the \zwlat.

\item
Let $a$ and $b$ be distinct vertices of $G$.  Suppose that $a$ and $b$ are both adjacent to $v$.  The vertices $a$ and $b$ are \newword{\oppsite\ neighbors}, respectively \newword{\nearby\ neighbors}, of $v$ if the edges $\simpoab va$ and $\simpoab vb$ are parallel, respectively perpendicular.

\item
The vertex $v$ is a \newword{\crnr} of $G$ if $v$ is adjacent to precisely two other vertices of $G$ and the two vertices adjacent to $v$ are \nearby\ neighbors.

\item
The \latgr\ $G$ is \newword{\twoless} if the degree of every vertex is at most $2$.
\qef
\eenum
\edefnnb

Clearly a graph being \twoless\ is equivalent to the graph having components that are isolated vertices, path graphs and cycle graphs.

\defnnb\label{defnBFR}
Let $R$ be a \ralp, and let $c$ be a crossing of $R$.
\enum
\item
The \newword{\xstrand}, respectively \newword{\ystrand}, of $R$ at $c$ is the union of two \xedges, respectively two \yedges, of $R$ that have $c$ as an endpoint.

\begin{figure}[ht]
\centering\includegraphics[scale=0.4]{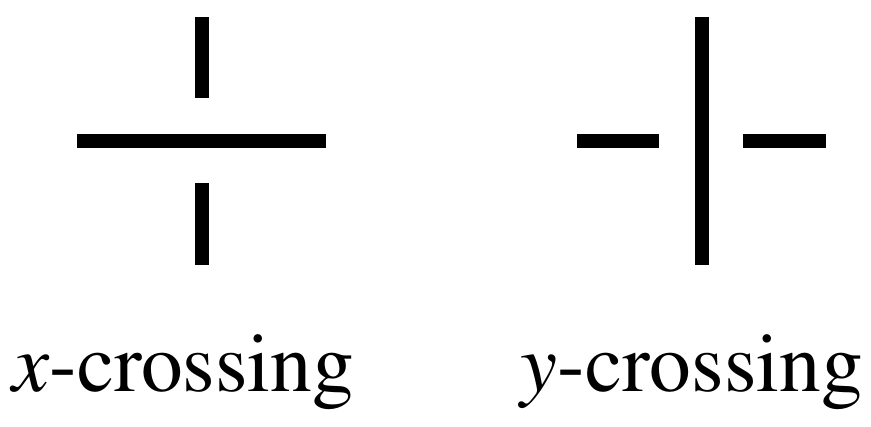}
\caption{}\label{vh_crossings}
\end{figure}

\item
The crossing $c$ is an \newword{\hcrossing}, respectively \newword{\vcrossing}, if the upper strand of the crossing is the \xstrand, respectively \ystrand, as seen in \figref{vh_crossings}.
\qef
\eenum
\edefnnb

\defn\label{defnBET}
Let $R$ be a \ralp.  The \newword{\crosgr} of $R$, denoted $\crg R$, is the lattice graph with a vertex at every crossing of $R$ and an edge between any two vertices that are adjacent in the \zwlat.
\edefn

\begin{figure}[ht]
  \centering
  \begin{minipage}[b]{0.45\textwidth}
    \centering
    \includegraphics[width=0.75\textwidth]{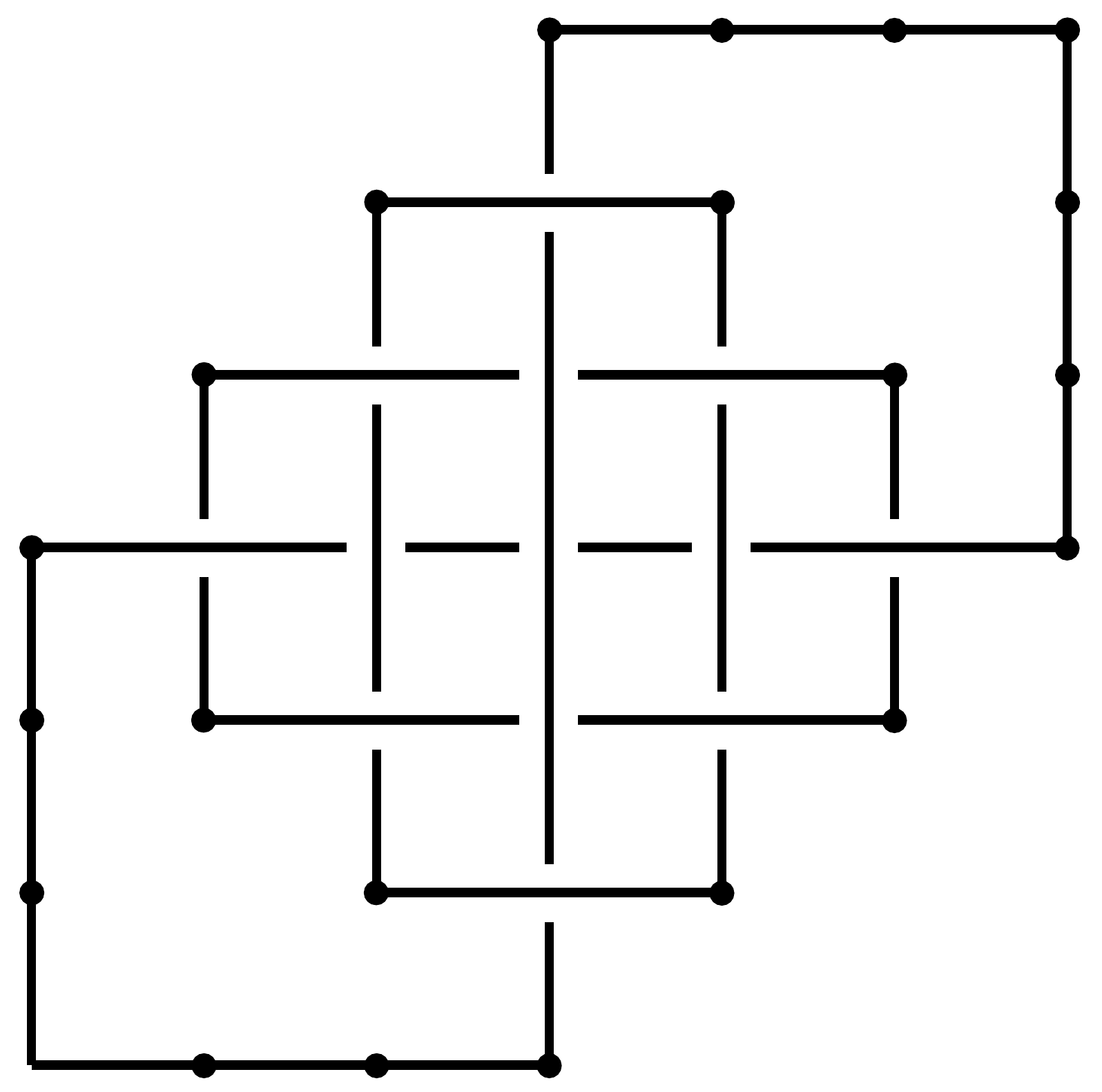}
    \caption{}\label{ralp_knot_pg1}
  \end{minipage}
  \hfill
  \begin{minipage}[b]{0.45\textwidth}
    \centering
    \includegraphics[width=0.6\textwidth]{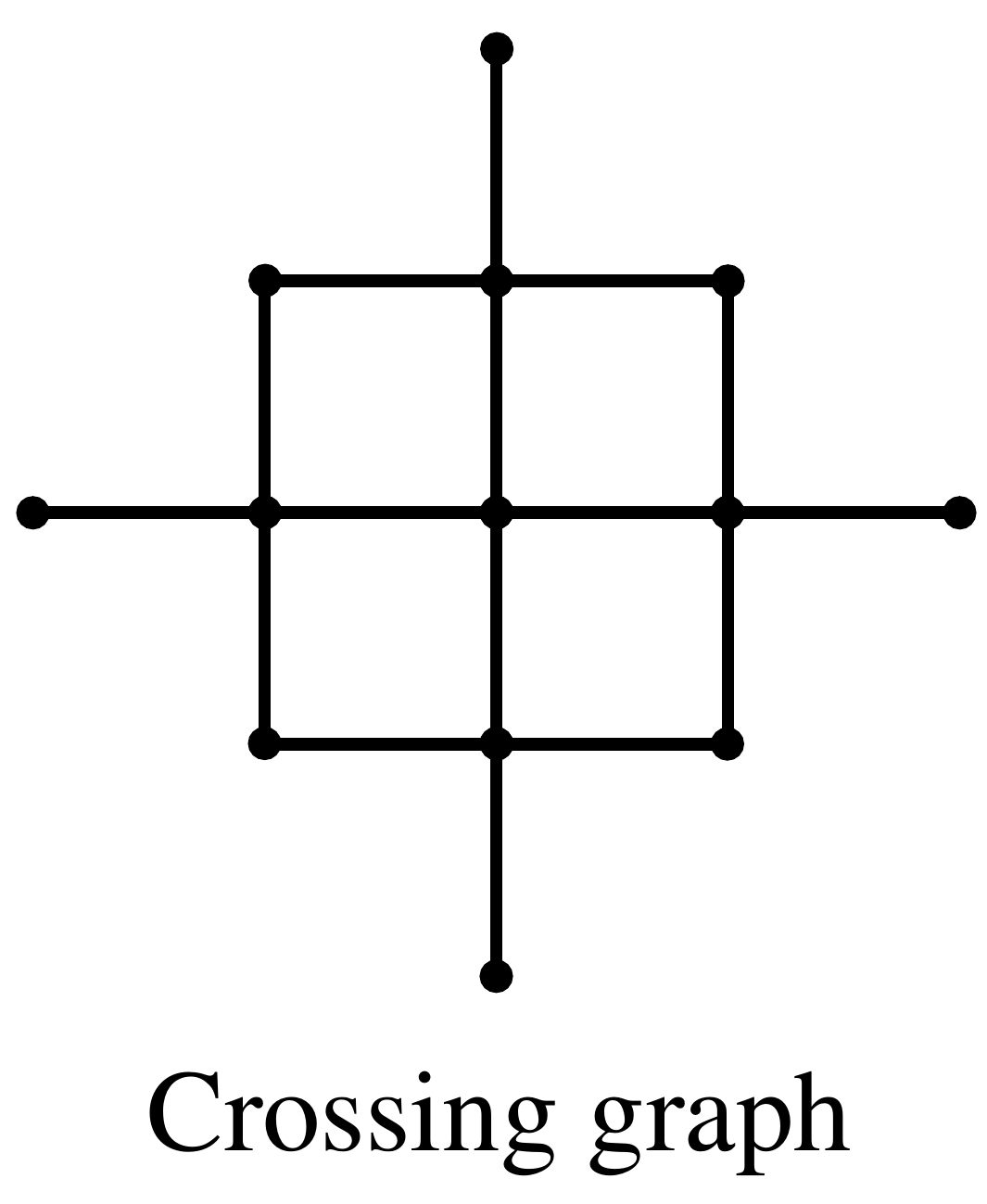}
    \caption{}\label{ralp_knot_pg2a}
  \end{minipage}
\end{figure}

For example, let $R$ be the \ralpl\ seen in \figref{ralp_knot_pg1}.  The \crosgr\ of $R$ is shown in \figref{ralp_knot_pg2a}.

\defn\label{defnBDE}
Let $R$ be a \ralp, and let $c$ be a crossing of $R$. The crossing $c$ is a \newword{\probcrs} if, thought of as a vertex in $\crg R$, it has a pair of \nearby\ neighbors, called \newword{\badnbrs} of $c$, which both have the opposite crossing type as $c$.
\edefn

For example, each of the crossings labeled $a$, $b$, $c$ and $d$ in \figref{ralp_knot_5_2a} are \probcrs, where the two labeled crossings next to any of these four vertices are its \badnbrs.

Observe that a \probcrs\ can have more than two \badnbrs.

The following remark is straightforward, and we omit the details.

\remknb\label{remkBDF}
Let $R$ be a \ralp.
\enum
\item\label{remkBDF1}
Let $d$ be a crossing of $R$.  Then $d$ is a not a \probcrs\ if and only if at least one pair of \oppsite\ neighbors of $d$ has the property that each of the two \oppsite\ neighbors either is not a crossing or is a crossing having the same crossing type as $d$.

\item\label{remkBDF2}
Suppose that four crossings of $R$ are the vertices of a unit square in the \zwlat.  These crossings form a \slpc\ if and only if each crossing of the four is a \probcrs\ and its two adjacent crossings among the four are among its \badnbrs.
\qex
\eenum
\eremknb

\lem\label{lemBDG}
Let $R$ be a \ralp.
\enum
\item\label{lemBDG1}
Suppose that two \probcrss\ in $R$ are adjacent.  Then either the two \probcrss\ have the same crossing type and neither is a \badnbr\ of the other, or they have the opposite crossing type and each is a \badnbr\ of the other. 

\item\label{lemBDG4}
Suppose that four crossings in $R$ are the vertices of a unit square.  If at least two of these four vertices is a \probcrs\ such that its two adjacent crossings among the four are among its \badnbrs, then all four of the crossings are \probcrss, and the four crossings form a \slpc.

\item\label{lemBDG5}
If $R$ has a \probcrs\ with more than two \badnbrs\ that are \probcrss, then $R$ has a \slpc.
\eenum
\elem

\demo
For Part~(\ref{lemBDG1}), let $c$ and $d$ be adjacent \probcrss.  If $c$ and $d$ have the same crossing type, then clearly neither is a \badnbr\ of the other, so suppose that $c$ and $d$ have the opposite crossing type.  Because $d$ is a \probcrs, then it must have at least two \badnbrs, all of which have the opposite crossing type as $d$.  If $c$ is not a \badnbr\ of $d$, then one of the \badnbrs\ of $d$, say $e$, is a \nearby\ neighbor of $d$ together with $c$, but that is a contradiction, because $e$ and $c$ both have the opposite crossing type as $d$, making $c$ a \badnbr\ of $d$.  Hence $c$ is a \badnbr\ of $d$.  A similar argument shows that $d$ is a \badnbr\ of $c$.

\begin{figure}[ht]
\centering\includegraphics[scale=0.4]{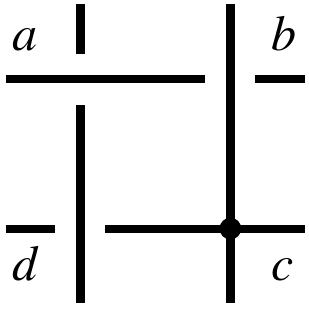}
\caption{}\label{ralp_basics_2_a}
\end{figure}

For Part~(\ref{lemBDG4}), let $a$, $b$, $c$ and $b$ be four crossings in $R$ that are the vertices of a unit square.  Suppose that at least two of these four vertices is a \probcrs\ such that its two adjacent crossings among the four are among its \badnbrs.  Without out loss of generality, suppose that $a$ is a \probcrs, and that $a$ is an \hcrossing\ such that its two adjacent crossings among the four are among its \badnbrs.  See \figref{ralp_basics_2_a}.  We know that $b$ and $d$ must be \vcrossings; the crossing shown by a dot is unspecified as of yet.  By hypothesis at least one of $b$, $c$ or $d$ is also a \probcrs\ such that its two adjacent crossings among the four are among its \badnbrs, and in any of these cases it is clear that $c$ must be an \hcrossing, which makes all four of the vertices \probcrss, and make the four crossings into a \slpc, by \remkrefp{remkBDF}2. 

\begin{figure}[ht]
\centering\includegraphics[scale=0.4]{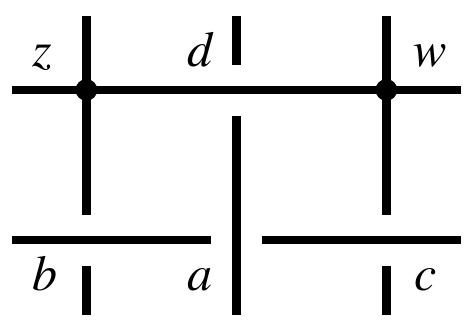}
\caption{}\label{ralp_basics_3}
\end{figure}

For Part~(\ref{lemBDG5}), suppose that $\pcrg R$ has a \probcrs\ $a$ with more than two \badnbrs\ that are \probcrss.  We consider the case where $a$ has exactly three \badnbrs\ that are \probcrss; the case where $a$ has four \badnbrs\ that are \probcrss\ is similar.  Without loss of generality, suppose that $a$ and its three \badnbrs\ that are \probcrss, denoted $b$, $c$ and $d$, are as seen in \figref{ralp_basics_3}, and that $a$ is a \vcrossing.  Hence $b$, $c$ and $d$ are \hcrossings; the crossings shown by dots are unspecified as of yet.  Clearly $z$ and $w$, as shown in the figure, are also crossings.

By Part~(\ref{lemBDG1}) of this lemma we know that $a$ is a \badnbr\ of $d$, and hence at least one of $z$ or $w$ is also a \badnbr\ of $d$.  Hence either the four vertices $a$, $d$, $z$ and $b$, or the four vertices $a$, $d$, $w$ and $c$, would have at least two vertices being a \probcrs\ such that its two adjacent crossings among the four are among its \badnbrs, and it follows from Part~(\ref{lemBDG4}) of this lemma that $R$ has a \slpc.
\edemo

The following definition makes sense by \lemrefp{lemBDG}1.

\defn\label{defnBDI}
Let $R$ be a \ralp.  The \newword{\pcrosgr} of $R$, denoted $\pcrg R$, is the the subgraph of $\crg R$ with a vertex at every \probcrs\ of $R$, and an edge between any two vertices that are \badnbrs\ of each other.
\edefn

\begin{figure}[ht]
\centering\includegraphics[scale=0.4]{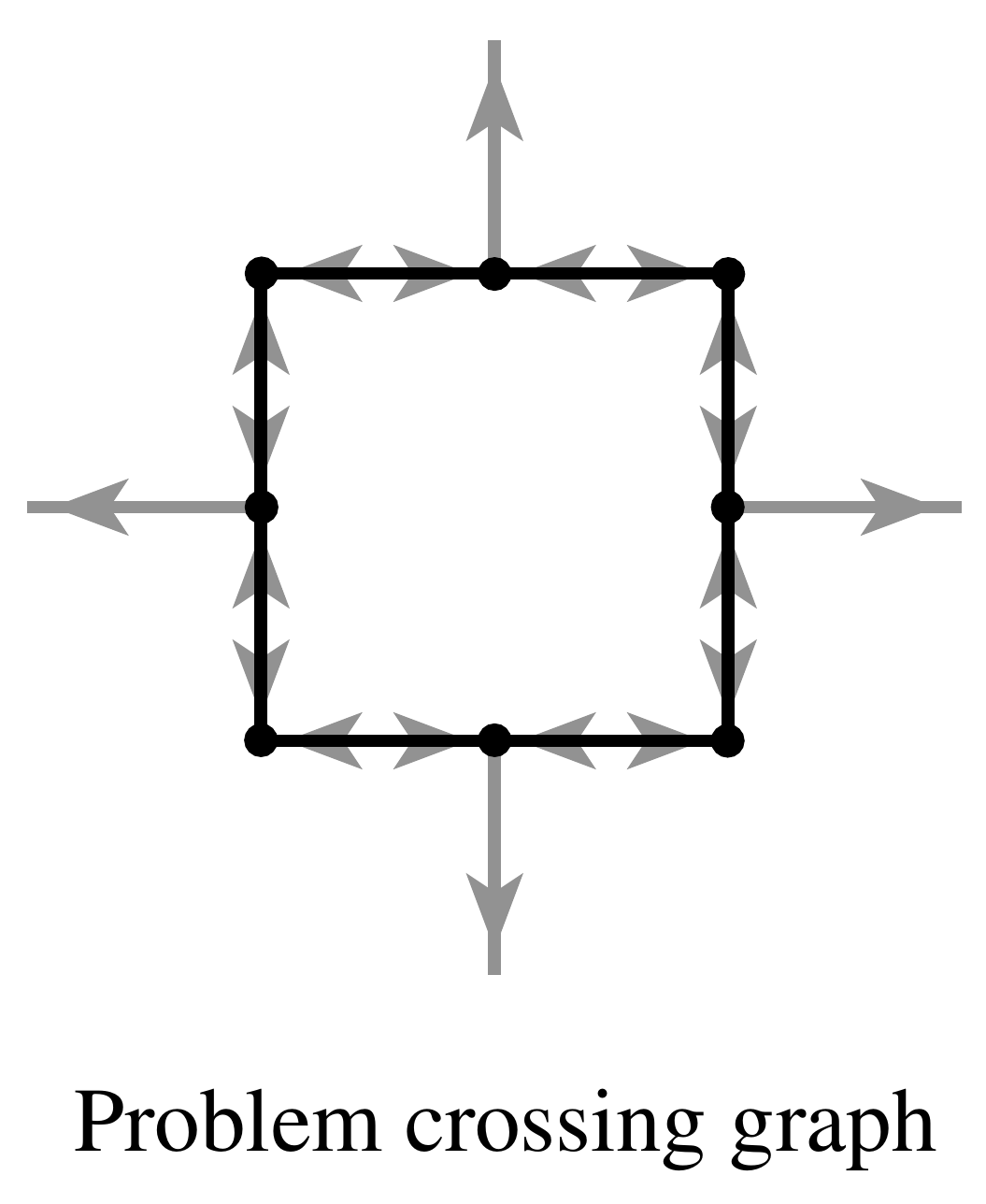}
\caption{}\label{ralp_knot_pg2b}
\end{figure}

If $R$ is the \ralpl\ seen in \figref{ralp_knot_pg1}, the \pcrosgr\ of $R$ is shown in \figref{ralp_knot_pg2b}, where the \pcrosgr\ itself consists of the edges and vertices in black, and where the gray arrows point from each \probcrs\ to its \badnbrs\ (some of which are in the \crosgr\ but not the \pcrosgr).

The following lemma is derived straightforwardly from \lemref{lemBDG}, combined with the fact that a \usf\ \latcyc\ must have a vertex that is not a \crnr; the proof is omitted.

\lem\label{lemBDH}
Let $R$ be a \ralp.  
\enum
\item\label{lemBDH1}
If two \probcrss\ of $R$ are adjacent in $\crg R$ but are not adjacent in $\pcrg R$, then neither is a \badnbr\ of the other, and they have the same crossing type.

\item\label{lemBDH2}
If $R$ does not have a \slpc, then $\pcrg R$ is \twoless\ and \usf, and every component of $\pcrg R$ that is a \latcyc\ has a vertex that is not a \crnr\ of $\pcrg R$.
\eenum
\elem


\section{Proof of the Theorem}
\label{secBASI}

To prove the non-trivial part of Theorem~\ref{thmBDM}, the idea is that we start with a \ralp, and we then modify it one crossing at a time, so that after each modification, the strands at that crossing do not intersect, and such that the crossings that were previously modified remain with their strands not intersecting; it is the latter that necessitates some care, which we accomplish by doing the modification at the crossings in a specific order.  After we complete all the modifications, we will end up with a \lsplskl\ that projects onto the \ralp.

Whereas a knot or link is an embedding of a simple closed curve or simple closed curves in $\rt$, it turns out that in the process of modification mentioned above we will have self-intersections of the knot or link, and hence we need a version of \lsplskls\ that is not necessarily embedded.  We start with the following preliminary.

\defnnb\label{defnBFS}
Let $C$ be a \tworeg\ graph.
\enum
\item
A function $\func fC{\rt}$ is a \newword{\latmap} if it is continuous, it maps every vertex of $C$ to a vertex of the \ztlat, and it maps every edge of $C$ onto a single edge of the \ztlat.

\item
Let $\func fC{\rt}$ be a \latmap.  An \newword{\xedge}, \newword{\yedge} or \newword{\zedge}, respectively, of $C$ with respect to $f$ is an edge of $C$ the image of which under $f$ is an \xedge, \yedge\ or \zedge, respectively, of the \ztlat.  An \newword{\xarc}, \newword{\yarc} or \newword{\zarc} of $C$ with respect to $f$ is a maximal arc of $C$ that is the union of \xedges, \yedges\ or \zedges, respectively.
\qef
\eenum
\edefnnb

Observe that if $C$ is a \tworeg\ graph and $\func fC{\rt}$ is a \latmap, then $C$ is the union of \xarcs, \yarcs, \zarcs, where such arcs intersect only in their endpoints.  To avoid special cases, a vertex in $C$ that is the endpoint of \xedges\ and/or \yedges, but not \zedges, can be considered to be a degenerate \zarc.

We now define our immersed version of \lsplskls, in relation to a given \ralp.

\defnnb\label{defnBFT}
Let $R$ be a \ralp, let $C$ be a \tworeg\ graph and let $\func fC{\rt}$ be a \latmap.  The \latmap\ $f$ is a \newword{\aplskl} for $R$ if the following four conditions hold.
\enum
\item[(a)]
$\zproj(f(C)) = R$.
 
\item[(b)]
The inverse image under $\zproj \rc f$ of any point in $R$ that is not a vertex in the \zwlat\ is a single point in $C$.

\item[(c)]
The inverse image under $\zproj \rc f$ of any vertex in $R$ that is not a crossing is a single \zarc\ in $C$ (possibly degenerate). 
 
\item[(d)]
The inverse image under $\zproj \rc f$ of any vertex in $R$ that is a crossing is two \zarcs\ in $C$ (each possibly degenerate), where one \zarc\ intersects only \xedges\ in $C$ and the other \zarc\ intersects only \yedges\ in $C$.
\qef
\eenum
\edefnnb

We note that a \ralp\ can be thought of an a \aplskl\ for itself.

\defnnb\label{defnBFQ}
Let $R$ be a \ralp.  Let $\func fC{\rt}$ be a \aplskl\ for $R$.  Let $v$ be a crossing of $R$.
\enum
\item
The \newword{\xzarc}, respectively \newword{\yzarc}, of $f$ at $v$, is the \zarc\ of $C$ in $(\zproj \rc f)^{-1}(v)$ that intersects only \xedges, respectively \yedges, in $C$.

\item
The \newword{\xzstick}, respectively \newword{\yzstick}, of $f$ at $v$, is the image under $f$ of the \xzarc, respectively \yzarc, of $f$ at $v$.  A \newword{\zstick} of $f$ at $v$ is either an \xzstick\ or a \yzstick. 

\item
The \newword{\xpstrand}, respectively \newword{\ypstrand}, of $f$ at $v$ is the image under $f$ of the \xzarc\ and the two \xedges\ of $C$ that the \xzarc\ intersects, respectively the image under $f$ of the \yzarc\ and the two \yedges\ of $C$ that the \yzarc\ intersects.

\item
The crossing $v$ is \newword{\resv\ with respect to $f$} if the following two conditions hold: (a) the \xzstick\ and \yzstick\ at $v$ are disjoint, and (b) if $v$ is an \hcrossing, respectively \vcrossing, then every point in the \xzstick\ at $v$ has larger, respectively smaller, \zvalue\ than every point in the \yzstick\ at $v$.
\qef
\eenum
\edefnnb

We note from \defnref{defnBFQ} that \xzsticks\ and \yzsticks\ are literal \zsticks\ in the \ztlat.  Additionally, we observe that the projection by $\zproj$ of the \xpstrand, respectively \ypstrand, of $f$ at $v$ is the \xstrand, respectively \ystrand, of $R$ at $v$.  Note also that the \xpstrand\ and \ypstrand\ of $R$ at $v$ each contains one of the two \zsticks\ over $v$.

\remk\label{remkBEX}
Let $R$ be a \ralp.
Let $\func fC{\rt}$ be a \aplskl\ for $R$.  Then the image of $f$ is a \lsplskl\ that projects onto $R$ if and only if each crossing of $R$ is \resv\ with respect to $f$.
\eremk

We now turn to the type of modification we use at the crossing of a \ralp.

We consider a \ralp\ to be in the \zwlat, and so prior to modification every crossing in a \ralp\ has \heig\ $0$.

\defnnb\label{defnBEQ}
Let $R$ be a \ralp.  Let $\func fC{\rt}$ be a \aplskl\ for $R$.  Let $v$ be a crossing of $R$.
\enum
\item
Let $p \in \zz$.  An \newword{\xlift{p}}\ at $v$ is a modification of $f$ defined as follows.  Let $s$ and $t$ be the two endpoints of the \xstrand\ of $R$ at $v$.  First, remove the \xzstick\ at $v$, and remove the \zstick\ at each of $s$ and $t$ that is connected to the \xstrand\ at $v$.  Second, move the images under $f$ of the two \xedges\ of $C$ in the \xpstrand\ at $v$ to the \heig\ $p$.  Finally, add in the necessary \zsticks\ to replace the three that were removed, in order to make the modified map be a \latmap\ (doing so may entail adding or removing edges from $C$).

\item
An \newword{\xlif}\ at $v$ is an \xlift{p}\ for some $p \in \zz$.

\item
A \newword{\ylift{p}}\ at $v$ for some $p \in \zz$, and a \newword{\ylif} at $v$, are defined similarly.

\item
A \newword{\lif} at $v$ is either an \xlif\ or a \ylif.

\item
A \lif\ at $v$ is \newword{\prlif} if after the \lif, the crossing $v$ is \resv\ with respect to $f$.\qef
\eenum
\edefnnb

\begin{figure}[ht]
  \centering
  \begin{minipage}[b]{0.45\textwidth}
    \centering
    \includegraphics[width=0.35\textwidth]{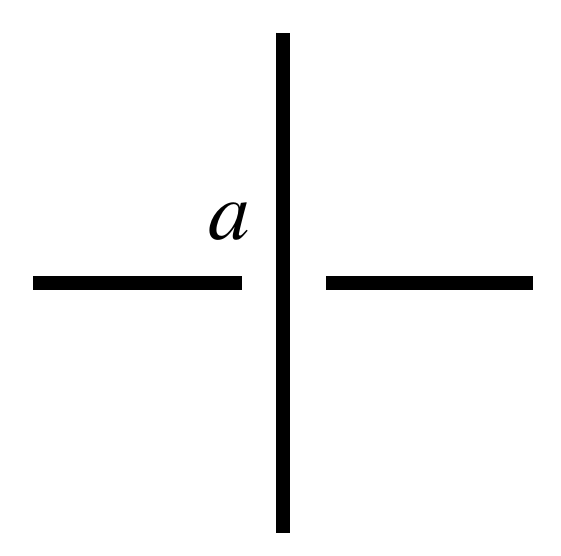}
    \caption{}\label{ralp_knot_no_lift_3a}
  \end{minipage}
  \hfill
  \begin{minipage}[b]{0.45\textwidth}
    \centering
    \includegraphics[width=0.35\textwidth]{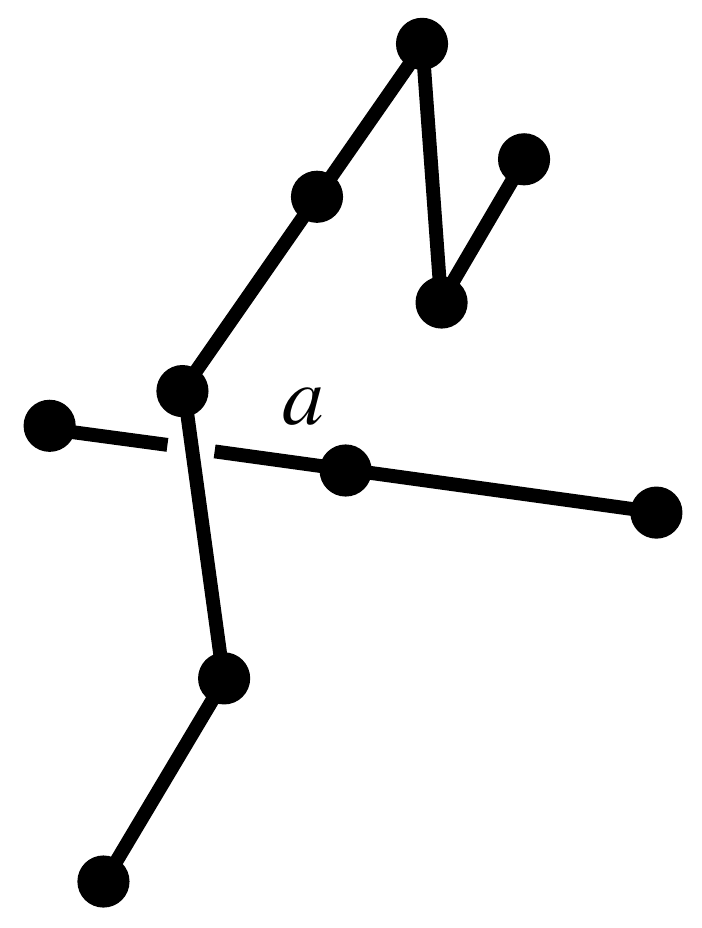}
    \caption{}\label{ralp_knot_no_lift_3b}
  \end{minipage}
\end{figure}

An example of a \pylift{1} is seen in \figrefs{ralp_knot_no_lift_3a}{ralp_knot_no_lift_3b}, where the former shows the crossing prior to the \lif, and the latter shows the crossing after the \lif.

Clearly, if $v$ is an \hcrossing\ of $R$, then doing an \xlift{p}\ or a \ylift{-p} for any sufficiently large $p \in \nn$ will be a \plif, and similarly for a \vcrossing.

We note that a \plif\ at $v$, while making $v$ be \resv\ or maintaining $v$ being \resv, might cause a neighboring crossing to go from being \resv\ to being not \resv.  For example, we see in \figref{ralp_basics_2_a} the crossings $a$ and $b$, which are an \hcrossing\ and a \vcrossing, respectively.  Suppose we did a \ylift{-1} at $a$, which is a \plif.  If we then did an \xlift{-1} at $b$, then $b$ would now be \resv, but $a$ would no longer be \resv.  Of course, if we did a \ylift{1} at $b$, that would make $b$ be \resv\ and would leave $a$ \resv; the problem would also be avoided if we had planned ahead and started with a \ylift{-2} at $a$, followed by an \xlift{-1} at $b$.  To avoid such problems, we use the following terminology.

\defnnb\label{defnBEI}
Let $R$ be a \ralp.  Let $\func fC{\rt}$ be a \aplskl\ for $R$.  Suppose that a \plif\ is done at a crossing $v$ of $R$.
\enum
\item\label{defnBEI1}
Let $w$ be another crossing of $R$.  The \lif\ at $v$ is \newword{\cmpt} with $w$ if the \lif\ at $v$ did not change $w$ from \resv\ to not \resv.  

\item\label{defnBEI2}
The \lif\ at $v$ is \newword{\bcmpt} if it is \cmpt\ with all other crossings.
\qef
\eenum
\edefnnb

Note that in \defnrefp{defnBEI}2, it does not matter what happens at a \nresv\ crossing other than $v$.  Note also that if a \lif\ at $v$ is an \xlif, respectively \ylif, then the only two crossings that might change from \resv\ to not \resv\ would be at the two endpoints of the \xstrand, respectively \ystrand, at $v$ (and only if those endpoints are crossings); in particular, if $w$ is a crossing of $R$ that is adjacent to $v$, and if the \lif\ at $v$ is perpendicular to $\simpoab vw$, then $w$ would not change from \resv\ to not \resv.

\begin{figure}[ht]
\centering\includegraphics[scale=0.38]{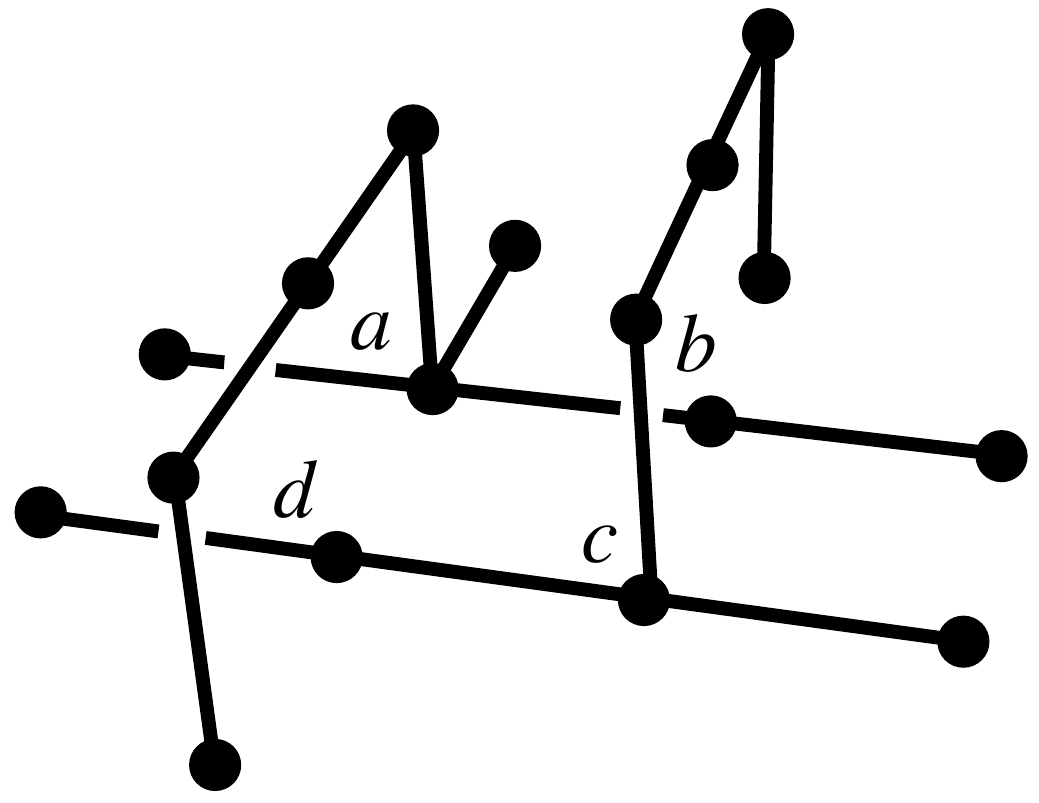}
\caption{}\label{ralp_knot_no_lift_2_a}
\end{figure}

We note that it is not always possible to do a \lif\ at a crossing that is \bcmpt.  The issue occurs at \probcrss.  For example, the crossing $a$ in \figref{ralp_basics_2_a} is a \probcrs, and suppose that a \ylift{+1}\ is done at each of $b$ and $d$, which makes these two crossings \resv, as seen in \figref{ralp_knot_no_lift_2_a}.  We then observe that a \plif\ at $a$ could be either an \xlift{p} for some $p \ge 2$, would change $b$ from \resv\ to not \resv, or a \ylift{-q} for some $q \ge 1$, which would change $d$ from \resv\ to not \resv.  Hence, no \plif\ is possible at $a$ that is \cmpt\ with both $b$ and $d$.

The above example is why the Proof of Theorem~\ref{thmBDM}, to which we now turn, is structured as it is.

\demoname{Proof of Theorem~\ref{thmBDM}}
If $R$ is has a \slpc, then the same argument used in regard to \figref{ralp_knot_5_2a} shows that $R$ is not the projection of a \lsplskl.

Now suppose that $R$ does not have a \slpc.  We can view $R$ as a \aplskl\ for itself.  We will do a \plif\ at one crossing of $R$ at a time, where each lift is \bcmpt.  After doing all the lifts, the resulting \aplskl\ will have all crossings \resv, and so it will be a embedding, and its image will be a \lsplskl\ that projects onto $R$.

We start with two preliminary observations about \lifs.  Let $v$ be a crossing of $R$.

Observation (1): Let $w$ be a crossing of $R$ that is adjacent to $v$.  If $v$ and $w$ have the same crossing type, then any \pxlif\ or \pylif\ at $v$ with sufficiently large \heig\ in absolute value is \cmpt\ with $w$.

To see why this observation is true, if $w$ is not \resv, there is nothing to prove, so suppose that $w$ is \resv.  Suppose further, without loss of generality, that $\simpoab vw$ is an \xedge.  If any \ylif\ is performed at $v$, then $w$ would not change from \resv\ to not \resv.  Hence, we need to consider only \xlifs\ at $v$.  There are two cases.

First, suppose that $v$ and $w$ are both \hcrossings.  Because $w$ is assumed to be \resv, then we note that the \xzstick\ at $w$ is higher than the \yzstick\ at $w$.  If an \xlift{n} is performed at $v$ where $n$ is larger than the highest point in the \xzstick\ at $w$, then such an \xlif\ at $v$ would not change $w$ from \resv\ to not \resv.  The case where $v$ and $w$ are \vcrossings\ is similar, except that the \xlif\ at $v$ has negative height.

Observation (2): Suppose $v$ is not a \probcrs.  Then there is a \lif\ at $v$ that is \bcmpt.

To see why this observation is true, we first note that by \remkrefp{remkBDF}1, at least one pair of \oppsite\ neighbors of $v$ has the property that each of these two \oppsite\ neighbors either is not a crossing or is a crossing having the same crossing type as $v$; let $a$ and $c$ be such \oppsite\ neighbors of $v$.  If $a$ and/or $c$ is not a crossing, or is a crossing that is not \resv, then there is nothing to be considered regarding that vertex, so assume that $a$ and $c$ are crossings that are \resv.  Observing that $a$, $v$ and $c$ all have the same crossing type.  Without loss of generality, we assume that the strand at the crossing $v$ that is in the direction of $a$ and $c$ is the upper strand at all three of $a$, $v$ and $c$.  A sufficiently high upward lift of the upper strand at $v$ will cause $v$ to be \resv, and will not cause $a$ or $c$ to change from \resv\ to not \resv.

We now return to doing \plifs\ at one crossing of $R$ at a time, starting with the \probcrss; if $R$ has no \probcrss, then skip this step.

By \lemrefp{lemBDH}2 we know that $\pcrg R$ is \twoless.  Hence, the components of $\pcrg R$ are isolated vertices, \latpaths\ and \latcycs. 

We note that if two \probcrs\ of $R$ are adjacent in $\crg R$ but are in different components of $\pcrg R$, or are in the same component of $\pcrg R$ but are not adjacent in $\pcrg R$, then by \lemrefp{lemBDH}1 we know that the two \probcrs\ have the same crossing type, and hence we can apply Observation~(1), which tells us that any \pxlif\ or \pylif\ with sufficiently large \heig\ in absolute value at each of these crossings is \cmpt\ with the other.  That tells us that we can do \plifs\ for each component of $\pcrg R$ separately without worrying about the impact on the other components of $\pcrg R$, and also that within a single component of $\pcrg R$, when we do a \plif\ at a vertex, we need only be concerned about being \cmpt\ with the adjacent vertices in $\pcrg R$, if there are any.

We proceed one component of $\pcrg R$ at a time, in any order. Let $M$ be a component of $\pcrg R$.

First, suppose $M$ is an isolated vertex (in $\pcrg R$, not in $\crg R$).  Let $v$ be the single vertex in $M$.  Then by a previous observation, there are \pxlifs\ and \pylifs\ at $v$ that are \cmpt\ with all other \probcrss; we do any such \lif\ at $v$.

Second, suppose $M$ is a \latpath.  Let $v_1, \ldots, v_n$ be the vertices of $M$ in order from one end of the arc to the other.  First, do a \plif\ at $v_1$ that is perpendicular to $\simpoab {v_1}{v_2}$.  Next, do a \plif\ at each of $v_2, \ldots, v_n$, in that order, such that the \lif\ at $v_i$ is perpendicular to $\simpoab {v_{i - 1}}{v_i}$ for each $i \in \{2, \ldots, n\}$.  For each $i \in \{2, \ldots, n\}$, we note that the perpedicularity implies that it is always possible to do such a \lif\ at $v_i$ that is \cmpt\ with $v_{i - 1}$; because $v_i$ is not adjacent in $\pcrg R$ to any of $v_1, \ldots, v_{i - 2}$, then the \lif\ at $v_i$ is \cmpt\ with $v_1, \ldots, v_{i - 1}$.  Hence, we can do \plifs\ at all vertices of $M$ that are \cmpt\ with all other \probcrss.

Third, suppose $M$ is a \latcyc.  By \lemrefp{lemBDH}2 we know that $M$ has a vertex that is not a \crnr.  Let $w_1, \ldots, w_m$ be the vertices of $M$ in order around the cycle, where $w_m$ is not a \crnr\ of $\pcrg R$.  We then proceed exactly as in the case where $M$ was a \latpath.  The only difference between the present case and the \latpath\ case is that in the present case, we need to ask whether the \lif\ at $w_m$ is \cmpt\ with $w_1$.  By construction we know that the \lif\ at $w_m$ is perpendicular to $\simpoab {w_{m - 1}}{v_m}$.  However, because $w_m$ is not a \crnr\ of $\pcrg R$, then we also know that the \lif\ at $w_m$ is perpendicular to $\simpoab {w_m}{v_1}$, and that means that the \lif\ at $w_m$ is \cmpt\ with $w_1$.  Again, we see that we can do \plifs\ at all vertices of $M$ that are \cmpt\ with all other \probcrss.

By doing the above to each of the components of $\pcrg R$, we have done lifts so that all the \probcrss\ are \resv.  Finally, we do a lift at one \nprobcrs\ at a time, which by Observation~(2) can always be done in a way that is \bcmpt.

As stated above, we have now found a \lsplskl\ that projects onto $R$.
\edemoname

\end{document}